\documentclass[12pt]{article}

\usepackage{amsmath,amssymb, fullpage}
\usepackage{graphicx}

\newcommand{\lecture}[4]{
   \pagestyle{myheadings}
   \thispagestyle{plain}
   \newpage
   \setcounter{page}{1}
   \noindent
}
\newtheorem{theorem}{Theorem}
\newtheorem{lemma}{Lemma}
\newtheorem{proposition}{Proposition}
\newtheorem{claim}{Claim}

\newcommand{\N}{\mathbb{N}}
\newcommand{\Z}{\mathbb{Z}}

\newcommand{\C}{\mathbb{C}}
\newcommand{\T}{\mathbb{T}}
\newcommand{\I}{\mathbb{I}}

\def\a{{\mathbf{\alpha}}}
\def\l{{\mathbf{\lambda}}}
\def\m{{\mathbf{\mu}}}
\def\n{{\mathbf{\nu}}}

\title{The computation of Kostka numbers and Littlewood-Richardson 
coefficients is \#P-complete}
\author{Hariharan Narayanan \\ University of chicago \\ 
		email:hari@cs.uchicago.edu}

\begin{document}
\maketitle
\begin{abstract}
Kostka numbers and Littlewood-Richardson coefficients play an essential 
role in the representation theory of the symmetric groups and the special linear 
groups. There has been a significant amount of interest in their computation
(\cite{barvinok}, \cite{rassart}, \cite{rassart2}, 
\cite{billey}, \cite{cochet}). 
The issue of their computational complexity has been a question of folklore,
but was asked explicitly by E. Rassart \cite{rassart}.
We prove that the computation of either quantity is 
\#P-complete. The reduction to computing  Kostka numbers, is from the \#P-complete problem \cite{kannan} 
of counting the number of $2 \times k$ contingency tables having given row and
column sums. 
The main ingredient in this 
reduction is a correspondence discovered by D. E. Knuth 
\cite{knuth}.  The reduction to the problem of computing Littlewood-Richardson 
coefficients is from that of computing Kostka numbers.
\end{abstract}

\section{Introduction}
Let $\N = \{1, 2, \dots \}$ be the set of positive integers and $\Z_{\geq 0} = \N \cup \{0\}$.
Let $\l = (\lambda_1, \dots, \lambda_s) \in \N^s$,
$\lambda_1 \geq \lambda_2 \geq \dots \geq \lambda_s \geq 1$,
$\m  = (\mu_1, \dots, \mu_t) \in \Z_{\geq 0}^t$,
$\n = (\nu_1, \dots, \nu_u) \in \Z_{\geq 0}^u$ and
$\a = (\alpha_1, \dots, \alpha_v) \in \N^v$, 
$\alpha_1 \geq \dots \geq \alpha_v \geq 1$.
The Kostka number $K_{\l \m}$ and the Littlewood-Richardson coefficient
$c^\n_{\l \a}$ play an essential role in the representation theory 
of the symmetric groups and the special linear groups. Their combinatorial definitions
can be found in Section~\ref{definitions}. 
The issue of their computational complexity has been a question of folklore.
Recently, in \cite{rassart}, E. Rassart asked whether there exist fast 
(polynomial time) algorithms to compute
Kostka numbers and Littlewood Richardson coefficients (Question 1, page 99).
We prove that these two quantities are \#P-complete 
(see Theorems \ref{th:k}, \ref{th:c}), and thus
answer his question in the negative under the
hypothesis that a \#P-complete quantity cannot be computed in polynomial time.

In \cite{barvinok}, Barvinok and Fomin show
how the set of all non-zero $K_{\l \m}$ for a given $\m$ can be produced
in time that is polynomial in the total size of the input and output.
They also give a probabilistic algorithm running in time, polynomial in the 
total size of input and output, that computes the set of all non-zero 
Littlewood-Richardson coefficients $c^\n_{\l \m}$ given $\l$ and 
$\m$. 
In \cite{cochet}, methods for the explicit computation of the 
Kostka numbers and Littlewood-Richardson coefficients using vector 
partition functions are discussed. 

Combinatorially, the Kostka number $K_{\l \m}$ is the number
of Young tableaux that have shape $\l$ and content $\m$ (\cite{fulton}, 
page 25). The Littlewood-Richardson coefficient $c^\n_{\l \a}$ is the number of 
LR skew tableaux on the 
shape $\l*\a$ with content $\n$ (this follows from Corollary 2, (v),
page 62 and Lemma 1, page 65 of \cite{fulton}).
Representation theoretically, $K_{\l \m}$ is the multiplicity of the 
weight $\m$ in the representation $V_\l$ of the lie algebra $sl_{r+1}(\mathbb{C})$ 
of the special linear group having 
highest weight $\l$ and $c^\n_{\l \a}$ is the multiplicity of 
$V_\n$ in the tensor product $V_\l \otimes_{\C} V_\a$.
They also appear in the representation theory of the symmetric groups (see
chapter 7, \cite{fulton}).

While there are formulae for  $K_{\l \m}$ and  $c^\n_{\l \a}$ 
due to Kostant and Steinberg respectively (\cite{cochet}, \cite{billey}), 
the number of terms is, in general,
exponential in the bit-length of the input.
These numbers have interesting properties such as,
for fixed $\l$, $\m$, $\a$, $\n$, $K_{N\l N\m}$ 
and $c^{N\n}_{N\l N\a}$ (\cite{rassart2}) are polynomials in $N$.
Whether $K_{\lambda\mu} > 0$ can be answered in polynomial time 
(see proposition~\ref{prop:k_in_p}), and so can 
the question of whether $c^\nu_{\lambda \alpha} > 0$, though the latter is  
a highly non-trivial fact estabilished by Ketan Mulmuley  
and Milind Sohoni \cite{ketan}, 
and uses the proof of the Saturation 
Conjecture by Knutson and Tao \cite{knutson}.
This fact plays an important role in the approach to the $P$ vs $NP$ question 
\cite{GCT} due to Ketan Mulmuley and Milind Sohoni.

Let $\mathbf{a}=(a_1, a_2) \in \Z_{\geq 0}, a_1 \geq a_2$ and 
$\mathbf{b} = (b_1, \dots b_k) \in \Z_{\geq 0}^k$.
We reduce the \#P-complete problem (\cite{kannan}) of finding the number 
$|\I(\mathbf{a}, \mathbf{b})|$ of $2 \times k$ contingency tables  
that have row sums $\mathbf{a} := (a_1, a_2)$ and
column sums $\mathbf{b} := (b_1, \dots, b_k)$, to that of finding 
some $K_{\lambda \mu}$. 
We reduce this in turn to finding some 
$c^\nu_{\l \a}$, where $\lambda, \mu, \alpha$ and $\nu$ can be computed
in time polynomial in the size of $(\mathbf{a}, \mathbf{b})$.
The main tool used in the reduction to finding Kostka numbers 
is a correspondence found by Donald E. Knuth (\cite{knuth}) between the set 
$\I(\mathbf{a}, \mathbf{b})$ of contingency tables and pairs of tableaux 
having contents
$\mathbf{a}$ and $\mathbf{b}$ respectively.


\section{Preliminaries and Notation}\label{definitions}

A {\bf Young diagram} (\cite{fulton}, page 1) is a collection of boxes, 
arranged in left justified rows, such that from top 
to bottom, the number of boxes in a row is monotonically (weakly) decreasing. 
The first two shapes in Figure $1$ are Young diagrams.
A {\bf filling} is a numbering of the boxes of a Young diagram with 
positive integers, that are
not necessarily distinct.
A {\bf Young tableau} or simply {\bf tableau} is a filling such that 
the entries are 
\begin{enumerate}
\item weakly increasing from left to right across each row, and
\item strictly increasing from top to bottom, down each column.
\end{enumerate}
$P$ and $Q$, in Figure $2$, are Young tableau.
A {\bf skew diagram} is the diagram obtained from removing a smaller Young
diagram out of a larger one. The third shape in Figure $1$ is a skew shape.
A {\bf skew tableau} is a filling of the boxes of a skew diagram with 
positive integers, non-decreasing in rows, and strictly increasing in columns
(see Figure $5$).
If the number of boxes in the
$i^{th}$ row of a tableau, for $1 \leq i \leq s$ is $\lambda_i$ and 
$\l:= (\lambda_1, \dots, \lambda_s)$, 
it is said to have {\bf shape} $\l$.
If the tableau houses $\mu_j$ copies of $j$ for $j \leq t$ and 
$\m := (\mu_1, \dots, \mu_t)$, it is said to have {\bf content} $\m$.
Thus, in figure $2$, $P$ and $Q$ have the same shape $(5, 2)$, but 
contents $(3, 2, 2)$ and $(4, 3)$ respectively.

\begin{figure}\label{fig:gamma}
\includegraphics[height=1.6in]{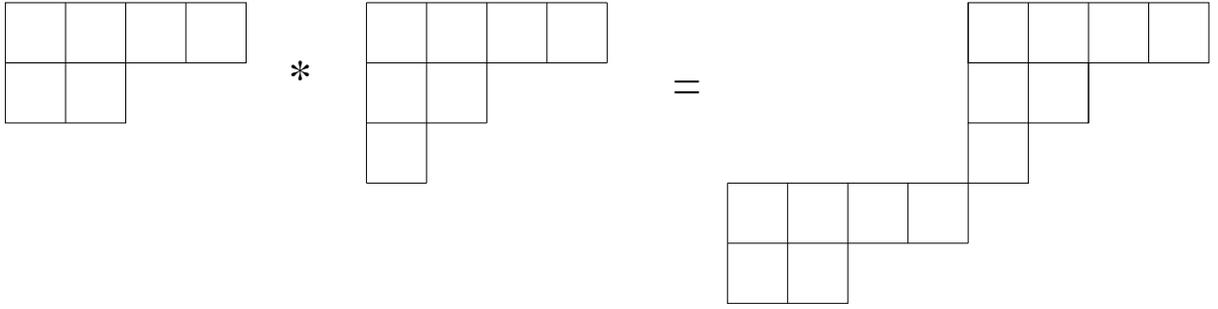}
\caption{Left to right, the shapes $\l, \a$ and the skew shape $\l*\a$.}
\end{figure}
We represent a Young diagram having $s$ rows and $\lambda_i$ boxes in the $i^{th}$ row
by $\l = (\lambda_1, \dots, \lambda_s)$. 
Given two shapes $\mathbf{\l}$ and $\mathbf{\a}$, 
$\mathbf{\l*\a}$ is defined to be the skew-shape 
obtained by attaching the lower left corner of $\a$ to
the upper right corner 
of $\l$ as in figure $1$ (see \cite{fulton}, page 60). 
${\bf size}(\lambda, \mu)$ denotes the number of bits used in the
description of this tuple of vectors.
For $\l := (\lambda_1, \dots, \lambda_s)$, 
let $|\l| = \sum_{i=1}^s \lambda_i$.
For vectors $\l$, $\m$, we say that $\l \unrhd \m$ if 
$|\l| = |\m|$ and $\forall i, \sum_{j \leq i} \l_j \geq \sum_{j \leq i} \m_j$. 
In addition, if $\l \neq \m$, we say $\l \rhd \m$. This ordering is called 
the dominance ordering.

We call a tableau {\bf Littlewood-Richardson} or {\bf LR}, 
if, when its entries are read right to left, 
top to bottom, at any moment, the number of copies of $i$ encountered is 
greater than or 
equal to the number of copies of $i+1$ encountered (\cite{fulton}, page 63). 
We denote the set of all (possibly skew) tableaux of shape $\l$ and content $\m$ by 
$\mathbf{\T(\l, \m)}$, and its subset consisting of all LR (possibly skew) tableaux
by LR$\T(\l, \m)$.
The {\bf Kostka number} $\mathbf{K_{\lambda\mu}}$ is the number of tableaux of 
shape $\l$ and content $\m$, i.e $|\T(\l, \m)|$ (\cite{fulton}, page 25).
The {\bf Littlewood-Richardson coefficient} $\mathbf{c^\n_{\l \a}}$ is
the number of LR skew tableaux of shape $\l*\a$ of content $\n$, i.e 
$|\text{LR}\T(\l*\a, \n)|$ (this follows from Corollary 2, (v),
page 62 and Lemma 1, page 65 of \cite{fulton}).

\section{The problems are in \#P}

Valiant defined the counting class \#P in his seminal paper 
\cite{valiant}, and proved that the computation of the permanent is 
\#P-complete.
The class \#P is the class of functions 
$f:\cup_{n \in \N}\{0, 1\}^n \rightarrow \Z_{\geq 0}$, for which there exists
a polynomial time turing machine $M$ and a polynomial $p$ such that
$(\forall n \in \N), (\forall x \in \{0, 1\}^n),  
f(x) = |\{y|y \in \{0, 1\}^{p(n)} \text{ and }  M \text{ accepts } (x, y)\}$.  

The tableau shapes $\l, \a$ and contents $\m, \n$ are described by vectors with integer coefficients. 
Therefore the number of boxes in a tableau might be exponential in the size of the input.
Thus it needs to be estabilished that even though an object in  
$\T(\l, \m)$ or $\text{LR}\T(\l*\a, \n)$ may be exponentially large, it has a polynomial length description,
and its membership (given this short description) can be verified in polynomial time.  

\begin{proposition}\label{prop:kp}
The problem of determining $K_{\l \m}$ is in \#P.
\end{proposition}
{\bf Proof:}\\
Let $\l, \m$ be defined as in the Introduction.
A tableau of shape $\l$ is fully specified given the number $\mu_{ij}$
of copies of $i$, that are present in the $j^{th}$ row, 
for all $i \leq s$ and $j \leq t$. 
Thus its description is of length polynomial 
in $\text{{\bf size}}(\l, \m)$, in fact it has length 
$O(\text{{\bf size}}(\l, \m)^2)$.   
Given the set $\{\mu_{ij}\}_{i \leq s, j \leq t}$, 
this corresponds to the tableau obtained, by writing, 
in the $i^{th}$ row, from left to right, 
$\mu_{i1}$ copies of $1$, $\mu_{i2}$ copies of $2$ and so on for each $i \leq s$. 
To verify whether this filling has the content $\m$, we need to verify, that
$\sum_i \mu_{ij} = \mu_j$ for each $j$, which takes $O(\text{{\bf size}}(\l, \m)^2)$ time.
To verify whether this filling has the shape $\l$, we need to verify, that
$\sum_j \mu_{ij} = \lambda_i$ for each $i$, 
which again takes $O(\text{{\bf size}}(\l, \m)^2)$ time.
To verify whether this filling is a tableau, we need to check 
that the entries in a column of this filling are strictly increasing from top to bottom. 
In other words, that, for each row $i \leq s-1$ and each $k \leq t-1$ 
$$\sum_{j \leq k}\mu_{ij} \geq \sum_{j \leq k+1} \mu_{i+1, j}.$$
This can be done in time $O(\text{{\bf size}}(\l, \m)^2)$ if we
maintained cumulative sums and calculated the next sum incrementally.
Therefore the task of computing $K_{\l \m}$ is in \#P.
\\$\Box$
\begin{proposition}\label{prop:k_in_p}
Given $\l$ and $\m$, whether or not $K_{\l \m} > 0$ can be answered in polynomial time.
\end{proposition}
{\bf Proof:}\\
Let $\l$ and $\m$ be as in proposition~\ref{prop:kp}. For any permutation
$\sigma$ of the set $\{1, \dots, t\}$, let $\sigma(\m)$ be the vector 
$(\mu_{\sigma(1)}, \dots, \mu_{\sigma(t)})$. 
It is a known fact that $K_{\l \m} = K_{\l \sigma(\m)}$ (see \cite{fulton}, page 26).
Let $\sigma$ be a permutation such that $\forall i \leq t-1$, 
$\mu_{\sigma(i)} \geq \mu_{\sigma(i+1)}$. For any $\check{\m}$, 
whose components are arranged in non-increasing order,
it is known that $K_{\l \check{\m}} > 0$ if and only if $\l \unrhd \check{\m}$ 
(see \cite{fulton}, page 26). 
Whether $\l \unrhd \sigma(\m)$ can be checked in time that is 
$O(\text{{\bf size}}(\l, \m))$.  
Thus, whether or not $K_{\l \m} > 0$ can be answered in time 
$O(\text{{\bf size}}(\l, \m) \ln(\text{{\bf size}}(\l, \m))$, 
which is the time it takes to find a permutation $\sigma$ that arranges the components of $\mu$ 
in non-increasing order.
\\$\Box$
\begin{proposition}\label{prop:cp}
The problem of computing $c^\n_{\l \a}$ is in \#P.
\end{proposition}
{\bf Proof:}\\
Let $\l, \a$ and $\n$ be defined as in the Introduction.
Given $S \in \text{LR}\T(\l*\a, \n)$, we shall describe it as follows.
Let the number of occurrences of $j$ in the $i^{th}$ row of $S$
be $\nu_{ij}$. 
$S$ is described by $(\l, \a)$ and $\{\nu_{ij}\}_{i \leq s+v,j \leq u}$.
Conversely, given $(\l, \a)$ and $\{\nu_{ij}\}_{i \leq s+v, j \leq u}$,
we verify that it has shape $\l*\a$ and content $\m$,
by checking that
\begin{enumerate}
\item $\forall i \leq v, \sum_j \nu_{ij} = \alpha_i$, 
\item $\forall i \geq v+1, \sum_j \nu_{ij} = \lambda_{i-v}$, and 
\item $\forall j, \sum_i \nu_{ij} = \nu_j$.
\end{enumerate}
This takes time $O(\text{size}(\l, \a, \n)^2)$. 
The description corresponds to the 
skew filling obtained, by writing in the shape $\l*\a$, left to right, 
in the $i^{th}$ row, 
$\nu_{i1}$ copies of $1$, $\nu_{i2}$ copies of $2$ and so on.
Such a skew filling would be a skew tableau if and only if the entries of each 
column strictly increased from the top to the bottom.
In other words,  
$$\sum_{j \leq k}\nu_{ij} \geq \sum_{j \leq k+1} \nu_{i+1, j},$$
for each row $i \leq s + v -1$ other than $s$ and each 
$k \leq u - 1$.
Thus, to verify from its description, that a given filling is an element of 
$\T(\l*\a, \n)$ takes polynomial time.
To verify that a specified filling is LR, we need to check, that, 
while scanning 
the entries of the tableau from left to right, top to bottom, at any instant, 
the number of copies of $j$ encountered is $\geq$ the number of copies of 
$j+1$ encountered.
In other words, that the 
following inequality holds for every row $k \leq v+s$ and every
$j \leq u-1$,
$$\sum_{i\leq k} \nu_{ij} \geq \sum_{i \leq k+1} \nu_{i, j+1}.$$
Thus, to verify that a skew tableau is
LR only takes time polynomial in the bit-length of its description.
Thus, the problem of determining $c^\n_{\l \a}$ is in \#P.
\\$\Box$

\section{Hardness Results}
\begin{lemma}\label{lemma:kp}
Given $\mathbf{a} = (a_1, a_2) \in \Z_{\geq 0}^2$, $a_1 \geq a_2$, and 
$\mathbf{b} = (b_1, \dots, b_k) \in \Z_{\geq 0}^k$, let $\l = (|a|, a_2)$ and
$\mu = (b_1, \dots, b_k, a_2)$. Then, 
$|\I(\mathbf{a}, \mathbf{b})| = K_{\lambda \mu}$.
\end{lemma}
{\bf Proof:}\\

\begin{figure}\label{fig:rsk}
\includegraphics[height=0.8in]{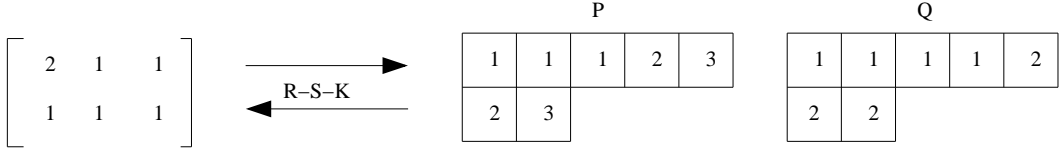}
\caption{An instance of the correspondence between 
$\I(\mathbf{a}, \mathbf{b})$ and 
$\cup_{\check{\l}}\T(\check{\l}, \mathbf{a}) \times \T(\check{\l}, \mathbf{b})$
for $\mathbf{a} = (4, 3)$, $\mathbf{b} = (3, 2, 2)$.}
\end{figure}
The R-S-K (Robinson-Schensted-Knuth) correspondence 
(\cite{knuth}, or \cite{fulton} pages 40-41) 
gives a bijection between $\I(\mathbf{a}, \mathbf{b})$, 
the set of $2 \times k$ contingency
tables with row sums $\mathbf{a}$ and column sums $\mathbf{b}$, 
and pairs of tableaux $(T_1, T_2)$
having a common shape but contents $\mathbf{a}$ and $\mathbf{b}$ respectively.
In other words, we have a bijection between  $\I(\mathbf{a}, \mathbf{b})$ 
and $\cup_{{\check{\l}}}\T(\check{\l}, \mathbf{a}) \times 
\T(\check{\l}, \mathbf{b})$.
A sample correspondence is shown in figure $2$.

\begin{figure}\label{fig:leaveP}
\includegraphics[height=0.8in]{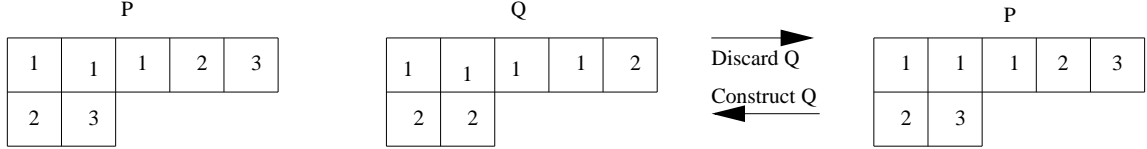}
\caption{An instance of the correspondence between
$\cup_{\check{\l}}\T(\check{\l}, \mathbf{a}) \times \T(\check{\l}, \mathbf{b})$
and $\cup_{\check{\l} \unrhd \mathbf{a}} \T(\check{\l}, \mathbf{b})$ for 
 $\mathbf{a} = (4, 3)$ and $\mathbf{b} = (3, 2, 2)$.}
\end{figure}

\begin{claim}
For every shape 
$\check{\l} = (\check{\lambda_1}, \check{\lambda_2})$, 
such that that $\check{\l} \unrhd \mathbf{a}$, 
there is exactly one tableau having shape $\check{\l}$ and content $\mathbf{a}$. 
For any other shape $\check{\l}$ there is no tableau having shape $\check{\l}$ 
and content $\mathbf{a}$.
\end{claim} 
Any tableau with content $\mathbf{a} = (a_1, a_2)$ can have at most two rows,
since the entries in a single column are all distinct.
Further, $1$ may be present only in the $1^{st}$ row, 
because top to bottom, each column has strictly increasing entries. 
Therefore for a tableau to have content $\mathbf{a}$, it is necessary that 
$|\l| = |\mathbf{a}|$ and there be
atleast $a_1$ boxes in the $1^{st}$ row. In other words, 
$$\T(\l, \mathbf{a}) \neq \phi \implies \check{\l} \unrhd \mathbf{a}.$$ 
Conversely, this condition together with $a_1 \geq a_2$, 
implies that  $a_1 \geq \check{\lambda_2}$. 
Therefore the filling in which the first 
$a_1$ boxes  of the top row contain $1$ and all others contain $2$ is 
a tableau (see $Q$ in Figure $3$). Since all the copies of $1$ must be in the first row and 
must be in a contiguous stretch including the leftmost box, this is the only
tableau in $\T(\l, \mathbf{a})$. Hence the claim is proved.
Thus there is a bijection between 
$\cup_{\check{\l}}\T(\check{\l}, \mathbf{a}) \times \T(\check{\l}, \mathbf{b})$
and the set of tableaux
of content $\mathbf{b}$ having some shape $\check{\l} \unrhd \mathbf{a}$. 
i.e, there is a bijection between 
$\cup_{\check{\l}}\T(\check{\l}, \mathbf{a}) \times \T(\check{\l}, \mathbf{b})$
and $\cup_{\check{\l} \unrhd \mathbf{a}} \T(\check{\l}, \mathbf{b})$. 
An example of this is provided in figure $3$.
Let us now consider the set 
$\cup_{\check{\l} \unrhd \mathbf{a}} \T(\check{\l}, \mathbf{b})$. 

\begin{claim}
Any tableau in  
$\cup_{\check{\l} \unrhd \mathbf{a}} \T(\check{\l}, \mathbf{b})$ 
can be extended to a tableau of the shape 
$\l = (a_1+a_2, a_2)$ by filling the boxes that are in $\l$ but not 
$\check{\l}$, with $k+1$.
This extension is a bijection between 
$\cup_{\check{\l} \unrhd \mathbf{a}} \T(\check{\l}, \mathbf{b})$ 
and $\T(\l, \m)$.
\end{claim}
If there is a tableau of shape $\check{\l}$ and content $\mathbf{a}$,  
$\check{\lambda_1} \leq |\mathbf{a}| = a_1 + a_2$, and 
$\check{\lambda_2} \leq a_2$.
$\check{\l} \unrhd \mathbf{a} \implies \check{\lambda_1} \geq a_2 = \lambda_2$.
Therefore no two of the boxes in $\l$ which are not in $\check{\l}$ 
belong to the same column. Those of these boxes, that are present in a given row,  
occupy a contiguous stretch that includes the rightmost box. Therefore
by filling them with $k+1$ we get a tableau in $\T(\l, \m)$.
Conversely, given a tableau $T$ in $\T(\l, \m)$, 
deleting all boxes of $T$ filled with $k+1$ 
gives a tableau in  
$\cup_{\check{\l} \unrhd \mathbf{a}} T(\check{\l}, \mathbf{b})$. 
These two maps are inverses of each other 
and hence provide a bijection
between  
$\cup_{\check{\l} \unrhd \mathbf{a}} T(\check{\l}, \mathbf{b})$ 
and $\T(\l, \m)$.
Hence the claim is proved.

An example of this correspondence has been illustrated in figure $4$.
\begin{figure}\label{fig:embed}
\includegraphics[height=0.8in]{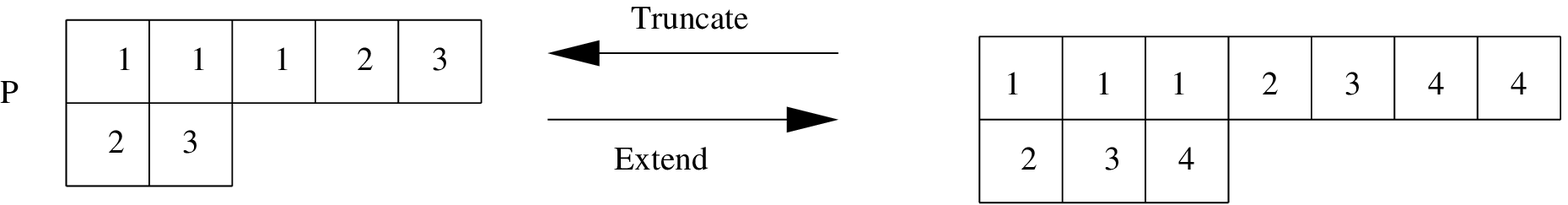}
\caption{An instance of the correspondence between 
$\cup_{\check{\l} \unrhd \mathbf{a}} \T(\check{\l}, \mathbf{b})$ and 
$\T(\l, \m)$, where $a = (4, 3), b=(3, 2, 2), \l = (7, 3)$ and 
$\m = (3, 2, 2, 3)$.}
\end{figure}
Therefore, $|\I(\mathbf{a}, \mathbf{b})|  
=|\cup_{\check{\l}}\T(\check{\l}, \mathbf{a}) \times \T(\check{\l}, \mathbf{b})|
=|\cup_{\check{\l} \unrhd \mathbf{a}} \T(\check{\l}, \mathbf{b})| = |\T(\l, \m)|
= K_{\l \m}$.
\\$\Box$

\begin{theorem}\label{th:k}
The problem of computing $K_{\l \m}$, even when $\l$ has only
$2$ rows, is \#P-complete.
\end{theorem}
{\bf Proof:}\\
Computing $K_{\l \m}$ is in \#P by proposition~\ref{prop:kp}.
Now the result follows from Lemma~\ref{lemma:kp} because  
the computation of $|\I(\mathbf{a}, \mathbf{b})|$ is known to be \#P-complete 
(\cite{kannan}).
\\$\Box$

\begin{lemma}\label{lemma:kc}
Given $\l = (\lambda_1, \lambda_2) \in \Z_{\geq 0}^2$, 
$\lambda_1 \geq \lambda_2$, and 
$\m = (\mu_1, \dots, \mu_\ell) \in \Z_{\geq 0}^\ell$, 
let $\a = (\alpha_1, \dots, \alpha_{\ell-1})$ where
$(\forall i) \alpha_i = \sum_{j>i} \mu_i$, 
and $\n = (\nu_1, \dots, \nu_{\ell})$,
where 
$\forall i \leq \ell - 1, \nu_i = \alpha_i + \mu_i$, and 
$\nu_\ell = \mu_\ell$.
Then $K_{\l \m} = c^{\n}_{\l \a}$.
\end{lemma}
{\bf Proof:}\\
\begin{figure}\label{fig:kc}
\includegraphics[height=1.3in]{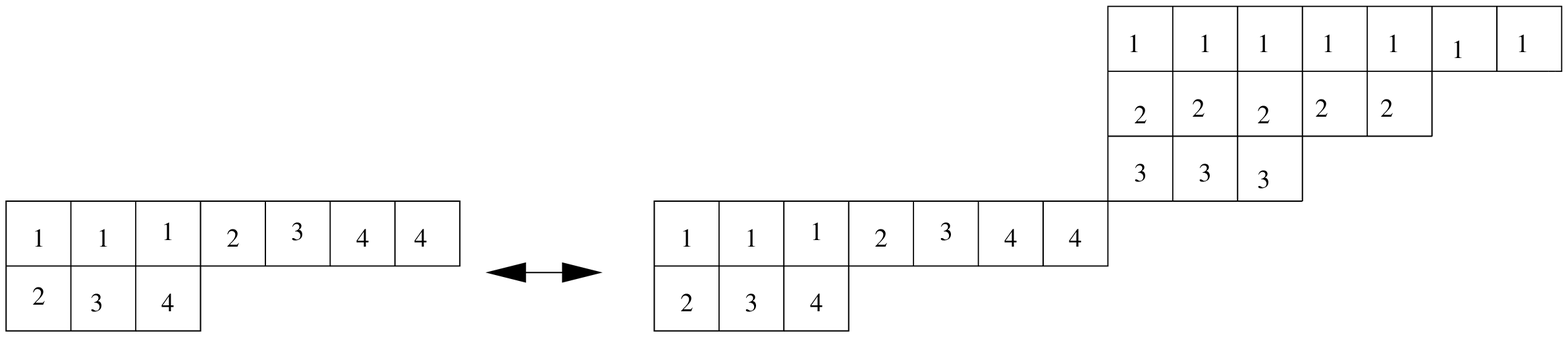}
\caption{An instance of the correspondence 
between $\T(\l, \m)$ and LR$\T(\l*\a, \n)$ for $\l = (7, 3)$ and $\m = (3, 2, 2, 3)$,
$\a = (7, 5, 3)$ and $\n = (10, 7, 5, 3)$.}
\end{figure}
$c^\nu_{\l \a}$ is, by definition, $|\text{LR}\T(\l*\a, \n)|$, 
which is the number of LR tableaux on the  skew shape
$\l*\a$ that have content $\n$.
The skew shape $\l*\a$ consists of a copy of $\l$ and a copy of $\a$, 
as in figures $1$ and $5$. For any skew tableau $S$ of shape $\l*\a$, 
we shall denote by
$S|_\a$, the restriction of $S$ to the copy of $\a$ and 
by $S|_\l$, the restriction of $S$ to the copy of $\l$.
Thus, $S|_\a$ is a tableau of shape $\a$ and $S|_\l$ is a tableau of shape $\l$. 

\begin{claim}
Let $S \in \text{LR}\T(\l*\a, \n)$. For $i \leq \ell-1$, the 
$i^{th}$ row of $S|_\a$ must consist entirely of copies of $i$.
\end{claim}
We proceed to show this by induction. 
Let $i = 1$.
The rightmost entry of the $1^{st}$ row of  $S|_\a$, is $1$,
by the LR condition.
The entries of each row are 
non-decreasing, left to right. Therefore each entry of the $1^{st}$ row must
be $1$.
Now assume for some $r$, $1 \leq r \leq \ell-1$, 
that the claim is true for all $i \leq r$.
If $r = \ell - 1$, we are done, so let $r < \ell - 1$. 
The entries of columns of  $S|_\a$ are strictly increasing, top to bottom.
By the induction hypothesis, the entry directly above a square in the 
$(r+1)^{th}$ row is $r$. Therefore, any element in the $(r+1)^{th}$ row must 
be $\geq r+1$. 
But the rightmost entry of the $(r+1)^{th}$ row cannot be anything $> r+1$,
because this would violate the LR condition. The entries of each row are 
non-decreasing, left to right. Therefore each entry of the $(r+1)^{th}$ 
row must
be $r+1$. Hence, by induction, the claim is proved.

Consequently, $S|_\l$ must have content $\n - \a = \m$. In other words,
$S|_\l \in \T(\l, \m)$.
Conversely, given any tableau $T \in \T(\l, \m)$, let $S(T)$ be the skew tableau
of shape $\lambda*\alpha$ in which $S(T)|_\l = T$ 
and the $i^{th}$ row of $S(T)|_\a$
consists entirely of copies of $i$. 
While scanning $S(T)$ right to left, top to 
bottom, the LR condition could not possibly be violated while on a square of 
$S(T)|_\a$.
By the time we begin scanning the rightmost box of the first row of 
$S(T)|_\l$, for any $i \leq \ell -1$, the number of copies of $i$ encountered 
is already $\mu_{i+1}$ more than the number of copies of $i+1$ encountered.
Therefore the LR condition could not possibly be violated on any square of 
$S(T)|_\l$ either. Therefore $S(T) \in \text{LR}\T(\l*\a, \n)$.
$S(T)|_\l = T$,
thus we have a bijection between $\text{LR}\T(\l*\a, \n)$, 
the set of LR skew tableaux of shape $\l*\a$ having 
content $\n$ and $\T(\l, \m)$, the set of tableaux of shape $\l$ having 
content $\m$.
Hence  $K_{\l \m} = |\T(\l, \m)|= |\text{LR}\T(\l*\a, \n)| = c^{\n}_{\l \a}$ as claimed. 
\\$\Box$

\begin{theorem}\label{th:c}
The problem of computing  $c^{\n}_{\l \a}$, even when $\l$ has only $2$ rows
is \#P-complete.
\end{theorem}
{\bf Proof:}\\
By Proposition~\ref{prop:cp}, computing $c^\n_{\l \a}$ is in \#P.
We have already proved in Theorem~\ref{th:k}, that the computation of 
$K_{\l \m}$ is \#P complete. The result now follows from Lemma~\ref{lemma:kc}.
\\$\Box$

\section{Conclusion}

We proved that the computation of Kostka numbers and Littlewood-Richardson 
coefficients is \#P complete. The reduction to computing Kostka numbers 
was from the \#P complete problem \cite{kannan} of computing the number of contingency tables
having given row and column sums. The problem of computing Kostka numbers
was then reduced to that of computing Littlewood-Richardson coefficients.
It is not known whether there exist Fully Polynomial Randomized 
Approximation Schemes (FPRAS) to compute any of these quantities. It would 
be of interest to find FPRAS to compute Kostka Numbers and 
Littlewood-Richardson coefficients. These are conjectured to exist by K. 
Mulmuley and M. Sohoni, \cite{ketan} in their approach to the P 
vs NP question.

\section{Acknowledgements}
I wish express my gratitude to Ketan Mulmuley for suggesting the topic of this paper and
for many valuable discussions. Many thanks are due to Ravi Kannan for informing me about
\cite{kannan}.
I also sincerely thank Etienne Rassart, Rahul Santhanam, and  L\'{a}szl\'{o} Babai for helpful comments.

\end{document}